\documentclass[reqno]{amsart}
\usepackage{amsfonts,graphicx,float}

\restylefloat{figure}

\makeatletter                                                  %
\def\th@plain{\slshape}                                        %
\makeatother                                                   %

\newcommand{\oi}{[0,1]}

\newcommand{\tot}{\leftrightarrow}

\newcommand{\Zbb}{\mathbb{Z}}
\newcommand{\Qbb}{\mathbb{Q}}
\newcommand{\Rbb}{\mathbb{R}}

\newcommand{\Rp}{\mathbb{R}^+}

\newcommand{\CH}{\mathbf{CH}}
\newcommand{\MV}{\mathbf{MV}}
\newcommand{\one}{{\rm 1\mskip-4mu l}}
\newcommand{\llgroup}{$\ell$-group}
\newcommand{\llgroups}{$\ell$-groups}

\newcommand{\pfrak}{\mathfrak{p}}

\newcommand{\mfrak}{\mathfrak{m}}

\newcommand{\Ccal}{\mathcal{C}}
\newcommand{\Dcal}{\mathcal{D}}
\newcommand{\Scal}{\mathcal{S}}

\newcommand{\Luk}{\L ukasiewicz}

\newcommand{\To}{\Rightarrow}

\def\dotminus{\buildrel\textstyle\cdot\over\relbar}

\newcommand{\newword}[1]{\textsl{#1}}

\newcommand{\vect}[3]{#1_#2,\ldots ,#1_#3}

\newcommand{\abs}[1]{\lvert#1\rvert}

\DeclareMathSymbol{\upharpoonright}{\mathrel}{AMSa}{"16}
\let\restriction\upharpoonright
\DeclareMathSymbol{\nmid}{\mathrel}{AMSb}{"2D}

\DeclareMathOperator{\Spec}{Spec}
\DeclareMathOperator{\MaxSpec}{MaxSpec}
\DeclareMathOperator{\den}{den}

\DeclareMathOperator{\Free}{Free}
\DeclareMathOperator{\Aut}{Aut}

\DeclareMathOperator{\SL}{SL}

\DeclareMathOperator{\Fl}{F\ell}

\theoremstyle{plain}
\newtheorem{theorem}{Theorem}[section]
\newtheorem{lemma}[theorem]{Lemma}

\newtheorem{corollary}[theorem]{Corollary}

\theoremstyle{definition}
\newtheorem{definition}[theorem]{Definition}
\newtheorem{example}[theorem]{Example}

\newtheorem{problem}[theorem]{}

\begin{document}

\bibliographystyle{plain}

\sloppy

\title[Automorphisms of product logic]{The automorphism group\\
of falsum-free product logic}

\author[G. Panti]{Giovanni Panti}
\address{Department of Mathematics\\
University of Udine\\
via delle Scienze 208\\
33100 Udine, Italy}
\email{panti@dimi.uniud.it}

\begin{abstract}
A few things are known, and many are unknown, on the automorphism group of the free MV-algebra over $n-1$ generators. In this paper we show that this group appears as the stabilizer of $\one$ in the larger group of all automorphisms of the free cancellative hoop over $n$ generators. Both groups have a dual action on the same space, namely the $(n-1)$-dimensional cube. The larger group has a richer dynamics, at the expense of loosing the two key features of the McNaughton homeomorphisms: preservation of denominators of rational points, and preservation of the Lebesgue measure. We present 
here some basic results, some examples, and some problems.
\end{abstract}


\thanks{\emph{2000 Math.~Subj.~Class.}: 06D35; 37B05}

\maketitle

\section{Preliminaries}
Consider the usual product $t$-norm $\cdot$, but restrict its domain by discarding $0$ from the real unit interval. The corresponding residuum is $a\to b=\min(1,b/a)$, and $\bigl((0,1],\cdot,\to,1\bigr)$ is a residuated lattice. The usual machinery of many-valued logic applies, and we have a logical system akin to the product logic in~\cite[\S4.1]{hajek98}. However, there is a key point of contrast: since $0$ is not a truth-value, every $n$-variables term-definable function from $(0,1]^n$ to $(0,1]$ is continuous. As a matter of fact, this falsum-free product logic is the logic of cancellative hoops. 

A \newword{cancellative hoop} is an algebra $(A,+,\dotminus,0)$ such that $(A,+,0)$ is a commutative monoid and the following identities are satisfied:
\begin{align*}
x\dotminus x &= 0, \\
x+(y\dotminus x)&=y+(x\dotminus y), \\
(x\dotminus y)\dotminus z&= x\dotminus(y+z), \\
x\dotminus(x\dotminus y)&=y\dotminus(y\dotminus x), \\
(x+y)\dotminus y &= x.
\end{align*}

If one fixes $0<c<1$, then the exponential function to base $c$ is an order-reversing isomorphism
$$
\exp:\bigl(\Rp,+,\dotminus,0\bigr)\to
\bigl((0,1],\cdot,\to,1\bigr)
$$
between $\bigl((0,1],\cdot,\to,1\bigr)$ and the positive cone $\Rp=[0,\infty)$ of $\Rbb$ endowed with the ordinary sum and the truncated difference $a\dotminus b=\max(0,a-b)$; note that $\exp(a\dotminus b)=\exp(b)\to\exp(a)$.

The following facts are known and can be found in~\cite{BlokFerreirim}, \cite{cignolitor00}, and references therein.
\begin{enumerate}
\item $(\Rp,+,\dotminus,0\bigr)$ generates the variety of cancellative hoops; this is a precise formulation of our previous statement that falsum-free product logic is the logic of cancellative hoops. 
\item Every cancellative hoop is obtainable as the positive cone of a unique enveloping lattice-ordered abelian group (\llgroup, for short), by setting $a\dotminus b=0\lor(a-b)$. This correspondence is 1-1 up to isomorphism,
and is a categorical equivalence. Under restriction to the positive cones,
\llgroup\ homomorphisms correspond bijectively to hoop homomorphisms, and \llgroup\ ideals (i.e., kernels of \llgroup\ homomorphisms) correspond to hoop ideals. In the following we will identify cancellative hoops with positive cones of \llgroups\ without further ado.
\item The free $n$-generated cancellative hoop $\Free_n(\CH)$ is the hoop of all continuous positively-homogeneous piecewise-linear functions with integer coefficients from $(\Rp)^n$ to $\Rp$, under pointwise operations.
\end{enumerate}

In this paper we will study the automorphism group of $\Free_n(\CH)$ 
and its dual action on the space of maximal ideals. It will turn out that this latter space is homeomorphic to the $(n-1)$-dimensional cube $\oi^{n-1}$, and that the action is given by piecewise-fractional homeomorphisms with integer coefficients. The automorphism group of the $n-1$-generated free MV-algebra sits naturally inside $\Aut\Free_n(\CH)$ as the stabilizer of a specific element. Passing from this stabilizer to the full group one gains a much richer dynamics, at the expense of loosing the two key features of the McNaughton homeomorphisms: preservation of denominators of rational points, and preservation of the Lebesgue measure.

We assume familiarity with \llgroups, MV-algebras, Mundici's $\Gamma$ functor, and the representation of the finitely generated free \llgroups\ and MV-algebras in terms of positively-homogeneous piecewise-linear functions and McNaughton functions, respectively \cite{bkw}, \cite{beynon77}, \cite{andersonfei}, \cite{glass99}, \cite{mundicijfa}, \cite{CignoliOttavianoMundici00}. The present paper pursues the line of research in~\cite{pantigeneric}, \cite{pantidynamical}, \cite{pantiinvariant}, \cite{pantibernoulli}; some acquaintance with some of the above papers is a prerequisite as well.

\section{Another description of $\Free_n(\CH)$}

Let $n\ge1$, and let $G_n$ be the \llgroup\ enveloping the free MV-algebra on $n-1$ generators $\Free_{n-1}(\MV)$, i.e., the unique \llgroup\ such that $\Free_{n-1}(\MV)=\Gamma(G_n,\one)$, where $\one$ is the function on $\oi^{n-1}$ whose value is constantly $1$. The elements of $G_n$ are all McNaughton functions from $\oi^{n-1}$ to $\Rbb$.
Let $\vect x1{{n-1}}\in G_n$ be the projection functions $\oi^{n-1}\to\oi$, and let $x_n=\one-(x_1\lor\cdots\lor x_{n-1})\in G_n$ (if $n=1$, then $\{\vect x1{{n-1}}\}$ is empty, and $x_1=\one$).

\begin{theorem}
The\label{ref1} free cancellative hoop over $n\ge1$ generators $\Free_n(\CH)$ is isomorphic to the positive cone $G_n^+$ of $G_n$, with $\vect x1n$ as free generators.
\end{theorem}
\begin{proof}
The case $n=1$ is clear, since $\Free_0(\MV)=\{0,1\}$ and $G_0=\Zbb$. Let $n>1$, let $Q=(\Rp)^n$ be the positive octant of $\Rbb^n$, and let $P\subset Q$ be the polyhedral cone spanned positively by $\{a_1e_1+a_2e_2+\cdots+a_{n-1}e_{n-1}+e_n:\vect a1{{n-1}}\in\{0,1\}\}$,
where $\{\vect e1n\}$ is the standard basis of $\Rbb^n$.
Denote by $\Fl_n$ the free \llgroup\ over $n$ generators: its elements are all continuous positively-homogeneous piecewise-linear functions with integer coefficients from $\Rbb^n$ to $\Rbb$.
In more detail, a function $F:\Rbb^n\to\Rbb$ is in $\Fl_n$ iff it is continuous and there exists a finite complex $\Sigma$ of rational polyhedral cones whose set-theoretic union is $\Rbb^n$ and such that, for each cone $W\in\Sigma$, there exist $\vect a1n\in\Zbb$ satisfying $F(\vect\alpha1n)=\sum a_i\alpha_i$ on $W$. Such a cone complex is usually called a \newword{fan}~\cite{fulton93}, \cite{ewald96}. A fan is \newword{unimodular} if all its cones are of the form $\Rp u_1+\cdots+\Rp u_t$, where $\vect u1t$ belong to $\Zbb^n$ and are extendable to a $\Zbb$-basis of $\Zbb^n$. The \newword{support} $\abs{\Sigma}$ of the fan $\Sigma$ is the set-theoretic union of all elements of $\Sigma$.

The projection functions $Y_i(\vect\alpha1n)=\alpha_i$ are free generators for $\Fl_n$.
Let now $I$ and $J$ be the principal ideals of $\Fl_n$ whose elements are all functions which are $0$ in $Q$ and in $P$, respectively.
The description of $\Free_n(\CH)$ given in \S1(3) amounts to saying that $\Free_n(\CH)$ is the positive cone $(\Fl/I)^+$ of $\Fl/I$, with free generators $Y_1/I,\ldots,Y_n/I$.
Note that the element $F/I$ of $\Fl_n/I$ is identifiable with the restriction $F\restriction Q$ of $F$ to $Q$, and analogously for $F/J$ and $F\restriction P$; we shall tacitly use such identifications.

For every element $\rho$ of the symmetric group over $n-1$ letters, let $N_\rho$ be the matrix obtained from the $n\times n$ matrix
$$
\begin{pmatrix}
0 & 1 & 1 & \cdots & 1 & 1 \\
0 & 0 & 1 & \cdots & 1 & 1 \\
0 & 0 & 0 & \cdots & 1 & 1 \\
\vdots & \vdots & \vdots & \cdots & \vdots & \vdots \\
0 & 0 & 0 & \cdots & 0 & 1 \\
1 & 1 & 1 & \cdots & 1 & 1
\end{pmatrix}
$$
by permuting its first $n-1$ rows according to $\rho$, and let $W_\rho$ be the unimodular cone spanned positively by the columns of $N_\rho$. Analogously, let $M_\rho$ be the matrix obtained from
$$
\begin{pmatrix}
0 & 1 & 1 & \cdots & 1 & 1 \\
0 & 0 & 1 & \cdots & 1 & 1 \\
0 & 0 & 0 & \cdots & 1 & 1 \\
\vdots & \vdots & \vdots & \cdots & \vdots & \vdots \\
0 & 0 & 0 & \cdots & 0 & 1 \\
1 & 0 & 0 & \cdots & 0 & 0
\end{pmatrix}
$$
by permuting its first $n-1$ rows according to $\rho$, and let $R_\rho$ be the cone spanned by the columns of $M_\rho$.
As proved in~\cite[Theorem~4.1]{pantigeneric}, the set of all $(n-1)!$ $W_\rho$'s and their faces is a unimodular fan $\Delta$ whose support is $P$, and the set of all $R_\rho$'s and their faces is a unimodular fan $\Sigma$ whose support is $Q$. The fans $\Delta$ and $\Sigma$ are combinatorially isomorphic, and by mapping every $W_\rho$ to the corresponding $R_\rho$ in the obvious way we obtain a map
$$
\Phi:P=\abs{\Delta}\to Q=\abs{\Sigma},
$$
which is an $\ell$-equivalence~\cite[p.~120]{beynon75}. By~\cite[pp.~120-121]{beynon75}, and taking into consideration the categorical equivalence of \S1(2), $\Phi$ induces an isomorphism of cancellative hoops
\begin{align*}
\varphi:(\Fl/I)^+&\to(\Fl/J)^+, \\
F\restriction Q &\mapsto F\circ \Phi\restriction P.
\end{align*}
By direct inspection one easily sees that, for every $1\le i\le n-1$ and every $W_\rho\in\Delta$, we have $Y_i\circ \Phi\restriction W_\rho=Y_i\restriction W_\rho$; hence 
$Y_i\circ \Phi\restriction P=Y_i\restriction P$ and
$\varphi(Y_i/I)=Y_i\circ \Phi\restriction P=Y_i\restriction P=
Y_i/J$ for every $1\le i\le n-1$.
We claim that the identity
$$
\bigl[(Y_1\lor\cdots\lor Y_{n-1})+Y_n\bigr]
\circ \Phi\restriction P=
Y_n\restriction P
$$
holds as well. Indeed, a typical point $u\in W_\rho$ is a column vector
$$
u=N_\rho
\begin{pmatrix}
\alpha_1 \\
\vdots \\
\alpha_n
\end{pmatrix}
$$
for certain $\vect\alpha1n\ge0$, and $Y_n(u)=\alpha_1+\cdots+\alpha_n$. On the other hand, we have
$$
\Phi(u)=M_\rho
\begin{pmatrix}
\alpha_1 \\
\vdots \\
\alpha_n
\end{pmatrix}
$$
and $\bigl[(Y_1\lor\cdots\lor Y_{n-1})+Y_n\bigr]\bigl(\Phi(u)\bigr)=(\alpha_2+\cdots+\alpha_n)+\alpha_1$, which settles our claim.
It follows that
\begin{align*}
\frac{Y_1\lor\cdots\lor Y_{n-1}}{J}+
\varphi\biggl(\frac{Y_n}{I}\biggr)&=
\varphi\biggl(\frac{Y_1\lor\cdots\lor Y_{n-1}}{I}\biggr)+
\varphi\biggl(\frac{Y_n}{I}\biggr) \\
&=\varphi\biggl(\frac{(Y_1\lor\cdots\lor Y_{n-1})+Y_n}{I}\biggr) \\
&=\frac{Y_n}{J},
\end{align*}
whence $\varphi(Y_n/I)=\bigl[Y_n-(Y_1\lor\cdots\lor Y_{n-1})\bigr]/J$.
Let $X_1=Y_1,\ldots,X_{n-1}=Y_{n-1},X_n=Y_n-(Y_1\lor\cdots\lor Y_{n-1})$;
we have shown that $\Free_n(\CH)$ is isomorphic to $(\Fl_n/J)^+$, with 
$X_1/J,\ldots,X_n/J$ as free generators. 

Observe now that the affine plane $\pi=\{Y_n=1\}$ cuts $P$ along a cross-section $\pi\cap P$ which is an $(n-1)$-dimensional cube. Since the elements of $(\Fl_n/J)^+$ are positively homogeneous, and the group and lattice operations act pointwise, the restriction $F/J\mapsto F\restriction(\pi\cap P)$ is an immersion of $(\Fl_n/J)^+$ into $G_n^+$; we 
write $f$ for $F\restriction(\pi\cap P)$.
The range of the immersion is freely generated by the images of $X_1/J,\ldots,X_n/J$, namely $\vect x1n$, and coincides with $G_n^+$ since, as it is well known, the latter is generated by $\vect x1{{n-1}},\one=x_n+(x_1\lor\cdots\lor x_{n-1})$.
\end{proof}

We have therefore three ways of looking at $\Free_n(\CH)$: as $(\Fl_n/I)^+$ with free generators $Y_1/I,\ldots,Y_n/I$, as $(\Fl_n/J)^+$ with free generators $X_1/J,\ldots,X_n/J$, or as $G_n^+$ with free generators $\vect x1n$; the first of these is essentially the one discussed in~\cite{cignolitor00}. 
In what follows, we always assume $n\ge2$ and, if not otherwise specified, we identify $\Free_n(\CH)$ with $G_n^+$ and $\Free_{n-1}(\MV)$ with $\Gamma(G_n,\one)$.

Let us remark that $\vect x1{{n-1}},\one$ generate $G_n^+$, but they are not free generators (except in the case $n=1$). Indeed, no strong unit $g$ can belong to a free generating set for $G_n^+$, since the map sending $g$ to $0$ and all other generators to $1$ cannot be extended to a homomorphism from $G_n^+$ to $\Rp$.

Theorem~\ref{ref1} makes clear the relationship between finitely generated free MV-algebras and free cancellative hoops. One obtains the free $n$-generated cancellative hoop by taking the positive cone of the \llgroup\ enveloping the free $(n-1)$-generated MV-algebra and forgetting about any distinguished strong unit. It might be slightly annoying that the $0$ function on $\oi^{n-1}$ means ``false'' in \Luk\ logic and ``true'' in falsum-free product logic, but things are promptly fixed by applying the $\neg$ involution to $\Free_{n-1}(\MV)$, so that $0$ means now ``true'' in both cases. It is reasonable to think of an element $f$ of $G_n^+$ to be ``sufficiently false'' if it is never true, i.e., has never value $0$ as a function from
$\oi^{n-1}$ to $\Rp$. Equivalently, $f$ is sufficiently false if it is a strong unit in $G_n$: for every $g\in G_n^+$ there exists an integer $m$ such that the conjunction of $f$ with itself $m$ times $f+\cdots+f$ is falser (i.e., greater) that $g$. One then goes back from $\Free_n(\CH)$ to $\Free_{n-1}(\MV)$ by deciding that a certain sufficiently false $f$ is actually the falsest proposition. The choice for $f$ is large, but not arbitrary. Indeed, $\Gamma(G_n,f)$ and $\Gamma(G_n,g)$ are isomorphic as MV-algebras iff there exists an automorphism of $G_n^+$ that maps $f$ to $g$. Therefore, $\Gamma(G_n,f)$ is isomorphic to $\Free_{n-1}(\MV)$ iff $f$ is in the orbit of $\one$ under the action of the automorphism group of $\Free_n(\CH)$. We shall see in~\S\ref{ref11} that this orbit is countably infinite, even though $\Aut\Free_n(\CH)$ does not act transitively on the set of strong units.

\section{Dual maps}

The \newword{spectral space}, or \newword{dual space}, of an \llgroup\ $G$ is the set of all the \newword{prime ideals} of $G$, i.e., the
kernels of nontrivial homomorphisms from $G$ to a totally-ordered group. This set is a topological space under the Zariski topology, in which a basic open set is the set of all kernels that avoid some fixed finite set of elements of $G$. The same description of the dual space applies both to MV-algebras and to commutative hoops: see~\cite{pantigeneric} for a more general and detailed presentation. We are interested here in the \newword{maximal spectrum} of the above structures, i.e., in the subspace of all kernels of nontrivial homomorphisms from $G$ to $\Rbb$ (to $\oi$ or to $\Rp$ in the case of MV-algebras or cancellative hoops, respectively). Fixing a strong unit in an \llgroup\ does not affect the spectrum, which is also preserved by the categorical equivalences between \llgroups, MV-algebras, and cancellative hoops. In short, for every \llgroup\ with strong unit $(G,u)$, the spectra of $G$, of $\Gamma(G,u)$, and of $G^+$, are identifiable in the obvious way, and the same holds for the maximal spectra.

The key point in the use of spectral spaces is the functoriality of the construction: to every endomorphism $\sigma$ of ---say--- the \llgroup\ $G$ it corresponds the \newword{dual map} $\sigma^*:\Spec G\to \Spec G$ given by $\sigma^*(\pfrak)=\sigma^{-1}[\pfrak]$; the dual map is automatically continuous.
A word of clarification is in order here: by definition, the trivial kernel $G$ is excluded from the spectrum. Hence, we must discard those endomorphisms $\sigma$ whose image $\sigma[G]$ is contained in some $\pfrak\in\Spec G$, for then $\sigma^*(\pfrak)$ would be undefined. We leave to the reader the proof of the following simple fact.

\begin{lemma}
Let $\sigma$ be an endomorphism of the \llgroup\ $G$. The following statements are equivalent:
\begin{itemize}
\item[(a)] for no $\pfrak\in\Spec G$ is $\sigma[G]\subseteq\pfrak$;
\item[(b)] for no $\mfrak\in\MaxSpec G$ is $\sigma[G]\subseteq\mfrak$;
\item[(c)] there exists a strong unit $u$ of $G$ such that $\sigma(u)$ is a strong unit;
\item[(d)] for every strong unit $u$ of $G$, the image $\sigma(u)$ is a strong unit.
\end{itemize}
An analogous statement holds for cancellative hoops.
\end{lemma}

Endomorphisms satisfying the above conditions are called \newword{nontrivial} in~\cite[p.~65]{pantigeneric}; note that endomorphisms 
of MV-algebras are automatically nontrivial. Nontrivial endomorphisms obviously form a monoid, and their dual maps are well defined, both on the spectrum and on the maximal spectrum.
Of course, any automorphism is nontrivial.

Our goal in this paper is to start an analysis of the group of automorphisms of $\Free_n(\CH)$ (or, rather, the dual group acting on the maximal spectrum). By way of comparison, let us sketch the situation
for MV-algebras. The maximal spectrum of $\Free_{n-1}(\MV)$ is homeomorphic to the $(n-1)$-cube $\oi^{n-1}$ with the standard Euclidean topology, under the identification
\begin{equation}\label{eq1}
\oi^{n-1}\ni p\mapsto\{f\in\Free_{n-1}(\MV):f(p)=0\}\in\MaxSpec\Free_{n-1}(\MV).
\end{equation}
Let $\sigma$ be any endomorphism of $\Free_{n-1}(\MV)$, and let $f_i=\sigma(x_i)$, for $1\le i\le n-1$. Then, under the
identification~(\ref{eq1}), the map $S$ dual to $\sigma$ turns out to be $S(p)=\bigl(f_1(p),\ldots,f_{n-1}(p)\bigr)$, and the mapping $\sigma\mapsto S$ is a contravariant monoid embedding: $\sigma\circ\tau\mapsto T\circ S$. A key point here is the existence of a distinguished measure on $\oi^{n-1}$ ---namely the Lebesgue measure--- which is preserved under the duals of all automorphisms~\cite{mundici95}. See~\cite{dinolagripa}, \cite{pantibernoulli}, \cite{pantiinvariant} for various results on the automorphism group of $\Free_{n-1}(\MV)$.

We now turn to cancellative hoops.
As we noted before, the maximal spectrum of $\Free_n(\CH)$ is identifiable with that of $\Free_{n-1}(\MV)$, i.e., with the $(n-1)$-cube. The correspondence~(\ref{eq1}) holds verbatim, by just substituting $\Free_{n-1}(\MV)$ with $\Free_n(\CH)$.

\begin{definition}
Let\label{ref6} $\sigma$ be a nontrivial endomorphism of $\Free_n(\CH)$, let $f_i=\sigma(x_i)$, for $1\le i\le n$, and let $f_\sharp=f_n+(f_1\lor\cdots\lor f_{n-1})=\sigma(\one)$. Since $\sigma$ is nontrivial, $f_\sharp$ is a strong unit, and never takes the value $0$. Moreover, $0\le f_i\le f_\sharp$ for every $1\le i\le n-1$. Let $S_i(p)=f_i(p)/f_\sharp(p)$;
the $n-1$ functions $S_i:\oi^{n-1}\to\oi$ are well defined.
We denote by $S$ the selfmap of $\oi^{n-1}$ given by $S(p)=\bigl(S_1(p),\ldots,S_{n-1}(p)\bigr)$.
\end{definition}

\begin{theorem}
Assume\label{ref2} the hypothesis and the notation of Definition~\ref{ref6}. Then, under the above identification of $\MaxSpec\Free_n(\CH)$ with $\oi^{n-1}$, the dual map of $\sigma$ is $S$.
\end{theorem}
\begin{proof}
Adopt all the definitions in the proof of Theorem~\ref{ref1}. In that proof we showed that $\Free_n(\CH)$ is isomorphic to $(\Fl_n/J)^+$ with free generators $X_1/J,\ldots,X_n/J$.
Let $u\in P$; then $u$ can be given two sets of coordinates: the $X$-coordinates $\bigl(X_1(u),\ldots,X_n(u)\bigr)$, and the $Y$-coordinates $\bigl(Y_1(u),\ldots,Y_n(u)\bigr)$.
The $Y$-coordinates are the ``real world coordinates'', i.e., the coordinates of $u$ in terms of the standard basis of $\Rbb^n$. The $X$-coordinates are the images of $u$ under the free generators of $(\Fl_n/J)^+$. If we write $\alpha_i=Y_i(u)$ and $\beta_i=X_i(u)$, then the two sets of coordinates are related by
\begin{align*}
\alpha_i &=\beta_i \text{ for $1\le i\le n-1$,}\\
\alpha_n &=\beta_n+(\beta_1\lor\cdots\lor\beta_{n-1}), \\
\beta_n  &=\alpha_n-(\alpha_1\lor\cdots\lor\alpha_{n-1}). 
\end{align*}
The functions in $\Fl_n$ are positively homogeneous, and the maximal spectrum of $(\Fl_n/J)^+$ is identifiable with the set of rays $\Rp u=\{ru:r\ge0\}$, for $u$ a point of $P$. Let $F_1/J,\ldots,F_n/J\in(\Fl_n)^+$ be such that the restrictions of $\vect F1n$ to $\pi\cap P$ are $\vect f1n$. 
Lift $\sigma$ to an endomorphism of $(\Fl_n/J)^+$, again denoted by $\sigma$, by setting $X_i/J\mapsto F_i/J$.
We define a selfmap $\Scal$ of $P$ as follows: 
if $u\in P$, then $\Scal(u)$ is the unique point of $P$ whose $X$-coordinates are $\bigl(F_1(u),\ldots,F_n(u)\bigr)$.
It is immediate that, for every $F/J\in(\Fl_n/J)^+$, we have $\bigl(\sigma(F)\bigr)(u)=F\bigl(\Scal(u)\bigr)$, and it follows that the map dual to $\sigma$ acts on $\MaxSpec(\Fl_n/J)^+$ by sending the ray through $u$ to the ray through $\Scal(u)$. The rays are in 1-1 correspondence with the points of the cross-section $\pi\cap P=\oi^{n-1}$ in the obvious way. Summing up, $S$ acts on $\oi^{n-1}$ as follows: given $p\in\pi\cap P$, construct the ray $\Rp\Scal(p)$ and intersect it with~$\pi$. The intersection contains a single point, namely $S(p)$.
In terms of coordinates, the $X$-coordinates of $\Scal(p)$ are $\bigl(f_1(p),\ldots,f_{n-1}(p),f_n(p)\bigr)$, and its $Y$-coordinates are $\bigl(f_1(p),\ldots,f_{n-1}(p),f_\sharp(p)\bigr)$. Dividing by $f_\sharp(p)$, we project $\Scal(p)$ to $\pi\cap P$, and we obtain the point having $Y$-coordinates $(S_1(p),\ldots,S_{n-1}(p),1)$, in accordance with the statement of the Theorem.
\end{proof}

\begin{corollary}
Assume~\label{ref10} the hypothesis and the notation of Definition~\ref{ref6}. Then, for every $f\in\Free_n(\CH)$, we have $\sigma(f)=f_\sharp\cdot(f\circ S)$.
\end{corollary}
\begin{proof}
Let $F/J$ be the unique element of $(\Fl/J)^+$ such that $f=F\restriction(\pi\cap P)$. Let $\Scal$ and $p$ be as in the proof of Theorem~\ref{ref2}. We have seen in that proof that $\Scal(p)=f_\sharp(p)\cdot S(p)$. Since $F$ is positively homogeneous, we obtain
$\bigl(\sigma(f)\bigr)(p)=\bigl(\sigma(F)\bigr)(p)=F\bigl(\Scal(p)\bigr)=F\bigl(f_\sharp(p)\cdot S(p)\bigr)=f_\sharp(p)\cdot
F\bigl(S(p)\bigr)=f_\sharp(p)\cdot(f\circ S)(p)$.
\end{proof}

A \newword{$k$-cell} $C$ in $\Rbb^n$ is a compact convex polyhedron of affine dimension~$k$. A \newword{rational cellular complex} $\Ccal$ is a finite set of cells 
such that: (1) all vertices of all cells in $\Ccal$ have rational coordinates; (2) if $C\in\Ccal$ and $D$ is a face of $C$, then $D\in\Ccal$; (3) every two cells intersect in a common face.
The \newword{support} of $\Ccal$ is the set-theoretic union $\abs{\Ccal}$ of all cells in $\Ccal$.

Let now $\sigma,\vect f1n,f_\sharp$ be as in Definition~\ref{ref6}. We can always partition the $(n-1)$-cube in a rational cellular complex $\Ccal$ in such a way that on each $(n-1)$-cell $C\in\Ccal$ all functions $\vect f1{{n-1}},f_\sharp$ are affine linear. Let us say that on $C$, and for $t\in\{1,\ldots,n-1,\sharp\}$, we have
$$
f_t=a^t_1x_1+\cdots+a^t_{n-1}x_{n-1}+a^t_n.
$$
Then, in homogeneous coordinates, $S\restriction C$ is given by
\begin{equation}\label{eq2}
\begin{pmatrix}
\alpha_1\\
\vdots\\
\alpha_{n-1}\\
1
\end{pmatrix}
\mapsto
\begin{pmatrix}
a^1_1 & \cdots & a^1_n \\
\vdots & \cdots & \vdots \\
a^{n-1}_1 & \cdots & a^{n-1}_n \\
a^\sharp_1 & \cdots & a^\sharp_n 
\end{pmatrix}
\begin{pmatrix}
\alpha_1\\
\vdots\\
\alpha_{n-1}\\
1
\end{pmatrix}.
\end{equation}
If $\sigma$ fixes $\one$, then $\sigma$ restricts to an endomorphism of the MV-algebra $\Gamma(G_n,\one)=\Free_{n-1}(\MV)$. In this case the last row of the above matrix is $(0\cdots0\,1)$, in accordance with~\cite[Theorem~2.6]{dinolagripa}.

The possible dynamics of dual maps in falsum-free product logic is far richer that that in \Luk\ logic. The following example shows some of the possibilities.

\begin{example}
Let\label{ref3} $a,b$ be positive integers. Let $f_1=b(x_1\land x_2)$, $f_2=a\bigl((x_1\lor x_2)\dotminus(x_1\land x_2)\bigr)$, and let $\sigma$ be the endomorphism of $\Free_2(\CH)$ determined by $x_i\mapsto f_i$. Since the $0$-sets of $f_1$ and $f_2$ do not intersect, $f_\sharp=f_1+f_2$ is a strong unit, and $\sigma$ is nontrivial. By explicit computation, one easily sees that the map $S$ dual to $\sigma$ depends only on the ratio $q=a/b$, and has the explicit form
$$
S(x)=\begin{cases}
x\cdot\bigl((1-2q)x+q\bigr)^{-1}, & \text{if $0\le x\le 1/2$;} \\
(1-x)\cdot\bigl((1-2q)(1-x)+q\bigr)^{-1}, & \text{if $1/2 < x\le 1$.}
\end{cases}
$$
We see here that, in contrast with the case of \Luk\ logic~\cite[p.~66]{pantigeneric}, the map $S$ does not determine $\sigma$: multiplying $a$ and $b$ for the same positive integer we get distinct $\sigma$'s and the same $S$. We plot here the graphs of $S$ for $q=2/9$, $q=1$, and $q=9/2$, respectively.

\begin{figure}[H]
\begin{center}
\includegraphics[height=3cm,width=3cm]{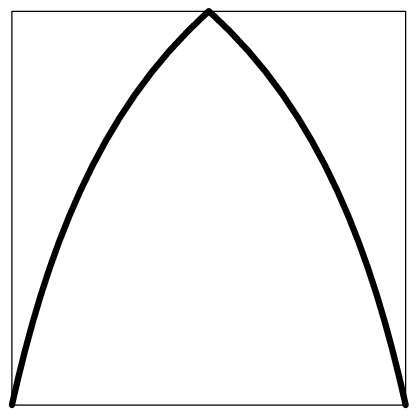}
\quad
\includegraphics[height=3cm,width=3cm]{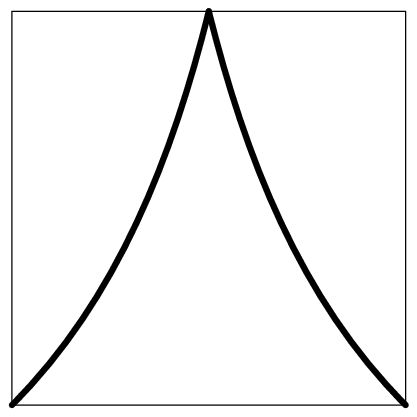}
\quad
\includegraphics[height=3cm,width=3cm]{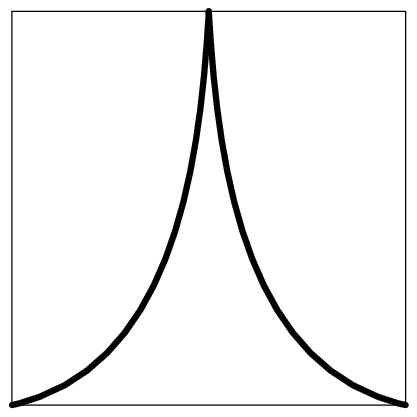}
\end{center}
\end{figure}

As shown in~\cite[Theorem~4.4]{pantigeneric}, we get quite different dynamics:
\begin{enumerate}
\item If $q<1$, then almost all (in the sense of Lebesgue measure) points have a dense $S$-orbit, and $S$ is ergodic with respect to a uniquely determined probability measure absolutely continuous w.r.t.~Lebesgue.
\item If $q=1$, then almost all points have a dense orbit, but $S$ does not preserve any probability measure absolutely continuous w.r.t.~Lebesgue.
\item If $q>1$, then almost all points are attracted to $0$.
\end{enumerate}
\end{example}

\section{The automorphism group}

We have seen in Example~\ref{ref3} that the map associating to a nontrivial endomorphism its dual map is not injective. Things go better in the case of automorphisms: we shall prove in Theorem~\ref{ref4} that if $S$ is the dual of an automorphism $\sigma$, then $S$ determines $\sigma$ uniquely.

\begin{definition}
Let\label{ref5} $S$ be an orientation-preserving homeomorphism of $\oi^{n-1}$. We call $S$ a \newword{piecewise $\SL_n\Zbb$-homeomorphism} if there exists a rational cellular complex $\Ccal$ whose support is $\oi^{n-1}$ and whose $(n-1)$-cells are $\vect C1k$, and there exist matrices $\vect A1k\in\SL_n\Zbb$ such that, for every $1\le h\le k$, $A_h$ expresses $S\restriction C_h$ in positively homogeneous coordinates (i.e., if $p=(\vect\alpha1{{n-1}})\in C_h$ and $S(p)=(\vect\beta1{{n-1}})$, then $A_h(\alpha_1\cdots\alpha_{n-1}\,1)^{tr}$ is positively proportional to $(\beta_1\cdots\beta_{n-1}\,1)^{tr}$).
\end{definition}

We call an automorphism of $\Free_n(\CH)$ \newword{orientation-preserving} or \newword{orientation-reversing} according whether its dual homeomorphism is orientation-preserving or orientation-reversing.
Of course the set of orientation-preserving automorphisms is a normal subgroup of $\Aut\Free_n(\CH)$ of index $2$, and every orientation-reversing automorphism is the composition of an orientation-preserving one with ---say--- the automorphism $x_1\tot x_2$.

\begin{theorem}
Let\label{ref4} $\sigma$ be an orientation-preserving automorphism of $\Free_n(\CH)$, with dual map $S$. Then $S$ is a piecewise $\SL_n\Zbb$-homeomorphism of $\oi^{n-1}$. Conversely, every piecewise $\SL_n\Zbb$-homeomorphism is the dual map of a unique orientation-preserving automorphism.
\end{theorem}
\begin{proof}
Taking into account the discussion after Corollary~\ref{ref10}, we need only prove that, if $C$ is an $(n-1)$-cell in $\oi^{n-1}$ in which $f_1=\sigma(x_1),\ldots,f_{n-1}=\sigma(x_{n-1}),f_\sharp=\sigma(\one)$ are all affine linear, then the matrix
$$
A=
\begin{pmatrix}
a^1_1 & \cdots & a^1_n \\
\vdots & \cdots & \vdots \\
a^{n-1}_1 & \cdots & a^{n-1}_n \\
a^\sharp_1 & \cdots & a^\sharp_n 
\end{pmatrix}
$$
in~\S3(\ref{eq2}) has determinant $1$. 
First of all, we extend uniquely $\sigma$ to an automorphism of $G_n$, again denoted by $\sigma$.
Let now $p=(\vect\alpha1{{n-1}})\in C$ be such that $\vect\alpha1{{n-1}},1$ are linearly independent over $\Qbb$, and let $\mfrak_p$ be the maximal ideal of $G_n$ whose elements are all McNaughton functions which are $0$ in $p$.
Evaluation at $p$ provides a canonical isomorphism from $G_n/\mfrak_p$ to the totally-ordered subgroup $H_p$ of $\Rbb$ generated by $\vect\alpha1{{n-1}},\alpha_n=1-(\alpha_1\lor\cdots\lor\alpha_{n-1})$ or, equivalently, by $\vect\alpha1{{n-1}},1$.
Let $q=S(p)=(\vect\beta1{{n-1}})$, and define analogously $H_q\simeq G_n/\mfrak_q$. Since $\sigma^{-1}[\mfrak_p]=\mfrak_q$, the map $f/\mfrak_q\mapsto\sigma(f)/\mfrak_p$ is an order isomorphism from $G/\mfrak_q$ to $G/\mfrak_p$. Denote by $\psi$ the corresponding order isomorphism from $H_q$ to $H_p$; by~\cite[Proposition~II.2.2]{kokorinkop}, $\psi$ must necessarily be of the form $\psi(\alpha)=r\alpha$, for a uniquely determined positive real number $r$. Since $\one/\mfrak_q\mapsto f_\sharp/\mfrak_p$, we have explicitly $r=f_\sharp(p)$.
It follows that $H_p$ coincides with the group $rH_q$ generated by the elements of the column vector
$$
r
\begin{pmatrix}
\beta_1\\
\vdots\\
\beta_{n-1}\\
1
\end{pmatrix}
=
\begin{pmatrix}
f_1(p)\\
\vdots\\
f_{n-1}(p)\\
f_\sharp(p)
\end{pmatrix}
=
A
\begin{pmatrix}
\alpha_1\\
\vdots\\
\alpha_{n-1}\\
1
\end{pmatrix}.
$$
Since $H_p=rH_q$ is isomorphic to $\Zbb^n$ as a group, $A$ must have either determinant $1$ or determinant $-1$, but the case $-1$ is excluded since $S$ is orientation-preserving.

Conversely, let $S$ be a piecewise $\SL_n\Zbb$-homeomorphism over a rational cellular complex $\Ccal$ as in Definition~\ref{ref5}. For every $1\le h\le k$, let $\Phi_h:\Rbb^n\to\Rbb^n$ be the nonsingular linear transformation whose associated matrix w.r.t.~the standard basis of $\Rbb^n$ is $A_h$. Since $A_h$ expresses $S\restriction C_h$ in positively homogeneous coordinates, $\Phi_h$ maps bijectively the cone $\Rbb^+C_h$ onto $\Rbb^+S[C_h]$. Moreover, if $p$ is a vertex common to $C_h$ and $C_t$, then $\Phi_h(p)=\Phi_t(p)$; this follows because $p$ has rational coordinates and $A_h,A_t\in\SL_n\Zbb$. Indeed, denoting by $u$ the \newword{primitive vector} along the ray $\Rbb^+p$ (i.e., the unique $u\in\Zbb^n\cap\Rbb^+p$ whose coordinates are relatively prime), then both $\Phi_h$ and $\Phi_t$ must map $u$ to the primitive vector along $\Rbb^+S(p)$. We conclude that the map $\Phi:P\to P$ ($P$ being the cone defined in the proof of Theorem~\ref{ref1}) defined by $\Phi(v)=\Phi_h(v)$ for $v\in\Rbb^+C_h$, is piecewise homogeneous linear with integer coefficients and hence, by~\cite[Corollary~1 to Theorem~3.1]{beynon77}, is induced by $n$ elements $\vect F1n$ of $\Fl_n$ as $\Phi(v)=\bigl(F_1(v),\ldots,F_n(v)\bigr)$.
By~\cite[Corollary~2 to Theorem~3.1]{beynon77} and the categorical equivalence between cancellative hoops and positive cones of \llgroups, there exists a unique automorphism of $(\Fl_n/J)^+$ (namely, the one defined by $X_i/J\mapsto F_i/J$) whose dual map on rays is $\Rbb^+v\mapsto\Rbb^+\Phi(v)$. Taking into account the isomorphism between $(\Fl_n/J)^+$ and $G_n^+$, and the correspondence between rays and $\MaxSpec(\Fl_n/J)^+$ in the proof of Theorem~\ref{ref2}, this concludes the proof of Theorem~\ref{ref4}.
\end{proof}

If $p=(\vect\alpha1{{n-1}})\in\oi^{n-1}\cap\Qbb^{n-1}$, the \newword{primitive homogeneous coordinates} of $p$ are the coordinates $(\vect a1{{n-1}},a_n)\in\Zbb^n$ of the primitive vector along the ray $\Rp(\vect\alpha1{{n-1}},1)$; the \newword{denominator} of $p$ is then $\den(p)=a_n$.

\begin{theorem}
Let~\label{ref13} $\sigma,S,\Ccal,\vect A1k$ be as in Theorem~\ref{ref4}. If $p$ is a point in the topological interior of some $(n-1)$-cell in $\Ccal$, then the Jacobian matrix $J(p)$ of $S$ at~$p$ is defined and its determinant has value $[f_\sharp(p)]^{-n}$.
The following statements are equivalent:
\begin{itemize}
\item[(a)] $\sigma(\one)=\one$ (i.e., $\sigma$ restricts to an automorphism of $\Free_{n-1}(\MV)$);
\item[(b)] for every $p\in\oi^{n-1}\cap\Qbb^{n-1}$, $\den\bigl(S(p)\bigr)=\den(p)$;
\item[(c)] for every vertex $p\in\Ccal$, $\den\bigl(S(p)\bigr)=\den(p)$;
\item[(d)] the last row of every $A_h$ is $(0\cdots0\,1)$;
\item[(e)] $S$ preserves the Lebesgue measure $\lambda$ on $\oi^{n-1}$ (i.e., $\lambda(T)=\lambda(S^{-1}T)$, for every measurable set $T$).
\end{itemize}
\end{theorem}
\begin{proof}
The statement about the Jacobian follows from~\cite[Proposition~2]{schweiger00}. Note that the set of points in which $S$ is not differentiable is contained in the union of the $(n-2)$-dimensional cells in $\Ccal$. The latter is a Lebesgue nullset, so we can safely  
write $\abs{J(p)}=[f_\sharp(p)]^{-n}$ throughout the $(n-1)$-cube
(recall that $f_\sharp$ is a strong unit, so it never takes value $0$).
The equivalence of (a) and (d) is the content of~\cite[Theorem~2.6]{dinolagripa}. (d) $\To$ (b) $\To$ (c) is clear;
we prove (c) $\To$ (d). Fix $1\le h\le k$, and choose $n$ vertices $\vect p1n$ of $C_h$ such that the matrix $B$ whose columns are the primitive homogeneous coordinates of $\vect p1n$ is nonsingular. Since $A_h$ has determinant $1$, the columns of $A_hB$ give the primitive homogeneous coordinates of $S(p_1),\ldots,S(p_n)$. By (c), we have the identity $(a_1\cdots a_{n-1}\, a_n)B=(0\cdots 0\,1)B$, where $(a_1\cdots a_{n-1}\, a_n)$ is the last row of $A_h$. Since $B$ is nonsingular, (d) follows. (a) $\To$ (e) is proved in~\cite[Theorem~3.4]{mundici95}. If (e) holds, then $\abs{J(p)}$ must be identically~$1$, so $f_\sharp^{-n}=[\sigma(\one)]^{-n}=\one$, and (a) holds as well.
\end{proof}

Let us say that a rational cellular complex $\Ccal$ supported in $\oi^{n-1}$ is a \newword{unimodular complex} if all $(n-1)$-cells in $\Ccal$ are simplexes and, for each such $(n-1)$-simplex~$C$, the primitive homogeneous coordinates of the vertices of $C$ constitute an $n\times n$ integer matrix whose determinant has absolute value $1$.
It is well known that every rational cellular complex can be refined to a unimodular one (see, e.g.,~\cite[Theorems~III.2.6 and VI.8.5]{ewald96}). It easily follows from
Theorem~\ref{ref4} that the most general automorphism of $\Free_n(\CH)$ is obtained by choosing a combinatorial isomorphism between
two unimodular complexes $\Ccal$ and $\Dcal$, both supported on $\oi^{n-1}$ (a \newword{combinatorial isomorphism} is a bijection between the two sets of vertices that preserves all incidence relations between simplexes).

\begin{example}
Take\label{ref7} $n=3$, and let $\Ccal$, $\Dcal$ be the following unimodular complexes:
\begin{figure}[H]
\begin{center}
\includegraphics[height=4.5cm,width=4.5cm]{figura-4}
\qquad
\includegraphics[height=4.5cm,width=4.5cm]{figura-5}
\end{center}
\end{figure}

\noindent The primitive homogeneous coordinates of the inner vertices are
\begin{align*}
&p_5=\begin{pmatrix}
1\\ 1\\ 3
\end{pmatrix},
&q_5=\begin{pmatrix}
2\\ 1\\ 4
\end{pmatrix},\\
&p_6=\begin{pmatrix}
1\\ 2\\ 5
\end{pmatrix},
&q_6=\begin{pmatrix}
1\\ 1\\ 3
\end{pmatrix}.
\end{align*}
We obtain a combinatorial isomorphism from $\Ccal$ to $\Dcal$ by mapping $p_t$ to $q_t$, for $1\le t\le 6$. If $C_1$ is ---say--- the simplex $\langle p_4,p_6,p_5\rangle$, then $A_1$ is uniquely determined by
$$
A_1\begin{pmatrix}
0 & 1 & 1 \\
1 & 2 & 1 \\
1 & 5 & 3
\end{pmatrix}
=
\begin{pmatrix}
0 & 1 & 2 \\
1 & 1 & 1 \\
1 & 3 & 4
\end{pmatrix}.
$$
Hence
$$
A_1=
\begin{pmatrix}
4 & 1 & -1 \\
2 & 2 & -1 \\
7 & 3 & -2
\end{pmatrix},
$$
and the restrictions of $f_1,f_2,f_3,f_\sharp$ to $C_1$ can be read from the rows of $A_1$:
\begin{align*}
f_1\restriction C_1&=4x_1+x_2-\one;\\
f_2\restriction C_1&=2x_1+2x_2-\one;\\
f_\sharp\restriction C_1&=7x_1+3x_2-2\one;\\
f_3\restriction C_1&=f_\sharp\restriction C_1-
(f_1\restriction C_1\lor f_2\restriction C_1).
\end{align*}
\end{example}

Note that the homeomorphism $S$ of Example~\ref{ref7} is not differentiable along the $1$-simplexes of $\Ccal$. Things are smoother in dimension~$1$, i.e., for $n=2$: first of all, a unimodular complex on $\oi$ is determined by a finite set of rational vertices
$$
\frac{0}{1}=\frac{a_0}{b_0}<
\frac{a_1}{b_1}<\cdots<
\frac{a_t}{b_t}=\frac{1}{1},
$$
such that $a_jb_{j+1}-a_{j+1}b_j=1$ for $0\le j\le t-1$. Two such complexes are combinatorially isomorphic iff they have the same number of vertices. If $S$ is the dual of some automorphism of $\Free_2(\CH)$, then $\abs{J(p)}=S'(p)$ for every $p$ in which $S$ is differentiable. By Theorem~\ref{ref13}, $S'$ equals either $f_\sharp^{-2}$ or $-f_\sharp^{-2}$, depending on whether $S$ is orientation-preserving or orientation-reversing. It follows that, in dimension~$1$, all dual homeomorphisms are of class $C^1$ (i.e., are differentiable everywhere with continuous first derivative).

\begin{example}
The\label{ref9} two sets of rational vertices
$$
\frac{0}{1}<
\frac{1}{2}<
\frac{2}{3}<
\frac{1}{1},
$$
and
$$
\frac{0}{1}<
\frac{1}{3}<
\frac{1}{2}<
\frac{1}{1},
$$
determine the map
$$
S(x)=\begin{cases}
x\cdot(x+1)^{-1}, & \text{if $0\le x\le 1/2$;} \\
(-x+1)\cdot(-5x+4)^{-1}, & \text{if $1/2< x\le 2/3$;} \\
(2x-1)\cdot x^{-1}, & \text{if $2/3 < x\le 1$;}
\end{cases}
$$
whose graph is
\begin{figure}[H]
\begin{center}
\includegraphics[height=4cm,width=4cm]{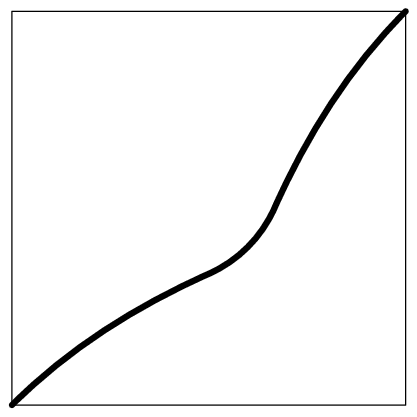}
\end{center}
\end{figure}

\noindent We see in this example that the duals of automorphisms of $\Free_2(\CH)$ may preserve no interesting probability measure.
As a matter of fact, basic ergodic theory shows that the only probability measures on $\oi$ 
preserved by $S$ are the affine combinations of the Dirac measures at $0$ and $1$.
\end{example}

\section{Open problems}

\begin{problem}
We\label{ref11} have seen that the automorphism group of the free MV-algebra over $n-1$ generators $\Aut\Free_{n-1}(\MV)$ is a (not normal) subgroup of $\Aut\Free_n(\CH)$, namely the stabilizer of $\one$. Our first problem is to describe the space of laterals, i.e., the orbit of $\one$ under $\Aut\Free_n(\CH)$. In equivalent terms, this means to characterize the set of strong units $g$ of $G_n$ such that $\Gamma(G_n,g)$ is isomorphic to $\Free_{n-1}(\MV)$. The orbit of $\one$ is countably infinite (it is already infinite under the action, e.g., of the single automorphism of Example~\ref{ref9}), but does not exhaust the set of strong units.
As an example, the strong unit $g=(x_1+\one)\land(-x_1+2\one)\in G_2^+$ is not in the orbit of $\one$, because $\Gamma(G_2,g)$ is not isomorphic to $\Free_1(\MV)$ (e.g., it has not a quotient isomorphic, as an MV-algebra, to $\{0,1/2,1\}$).
This should be compared with the case of real vector lattices, in which the automorphism group acts transitively on strong units; see~\cite[Lemma~4.2]{beynon75}.
\end{problem}

\begin{problem}
Given two rational points $p,q\in\oi^{n-1}\cap\Qbb^{n-1}$, is it always true that there exists $\sigma\in\Aut\Free_n(\CH)$ whose dual maps $p$ to $q$? Much harder: provided that $p$ and $q$ have the same denominator, can $\sigma$ be taken in $\Aut\Free_{n-1}(\MV)$?
\end{problem}

\begin{problem}
$\Aut\Free_{n-1}(\MV)$ is residually finite~\cite[p.~75]{dinolagripa}. What about $\Aut\Free_n(\CH)$? Is either group finitely generated?
\end{problem}

\begin{problem}
Let $n\ge3$, and let $g\in G_n^+$ be a strong unit in the orbit of $\one$, say $\sigma(g)=\one$. Then there exists a unique probability measure $\mu$ on $\oi^{n-1}$ which is null on underdimensioned $0$-sets (see~\cite[Definition~2.2]{pantiinvariant}) and is invariant under the duals of all automorphisms of $G_n^+$ that fix $g$. Namely, 
$\mu$ is the push-forward $S_*\lambda$ of the Lebesgue measure $\lambda$ via $S$: for every Borel set $A\subseteq\oi^{n-1}$, we have by definition $\mu(A)=(S_*\lambda)(A)=\lambda(S^{-1}A)$. These facts follow from the main result of~\cite{pantiinvariant}; see also the discussion in~\cite[\S2]{pantibernoulli}. We only show here the invariance of $S_*\lambda$: choose $\tau\in\Aut G_n^+$ such that $\tau(g)=g$. Then $\sigma\circ\tau\circ\sigma^{-1}$ fixes $\one$, and hence its dual $S^{-1}\circ T\circ S$ fixes the Lebesgue measure. Therefore $\lambda=(S^{-1}\circ T\circ S)_*\lambda=(S^{-1}_*\circ T_*\circ S_*)\lambda$, and $S_*\lambda=T_*(S_*\lambda)$, as desired. Now the question is: given a strong unit $g$ of $G_n^+$ which is \emph{not} in the orbit of $\one$, does it always exist a probability measure on $\oi^{n-1}$ which is invariant under the duals of all automorphisms that leave $g$ fixed? If so, under which conditions is such a measure unique?
\end{problem}

\begin{problem}
In~\cite{pantibernoulli} it is proved that for every odd $n$ there exists an automorphism of $\Free_{n-1}(\MV)$ whose dual is measure-theoretically isomorphic to a Bernoulli shift. What about even $n$?
\end{problem}

\end{document}